

\documentstyle{amsppt}    

\def\a{\kern+.6ex\lower.42ex\hbox{$\scriptstyle \iota$}\kern-1.20ex a}
\def\e{\kern+.5ex\lower.42ex\hbox{$\scriptstyle \iota$}\kern-1.10ex e}

\define\Ha{\Cal H}                       
\define\Cn{\Cal C^n}                     
\define\Pn{\Cal P^n}                     
\define\Ce{\Cal C}                       
\define\Hol{\Cal O}                      
\define\mult{\operatorname{mult}}        
\define\res{\operatorname{res}}          
\define\ord{\operatorname{ord}}          

\topmatter
\title The Noether exponent and Jacobi formula \endtitle
\author by Arkadiusz P\l{}oski \endauthor
\address
Department of Mathematics,\newline
Technical University,\newline
Al.~1000\,LPP\,7, 25--314~Kielce,\newline
Poland
\endaddress
\abstract
For any polynomial mapping $F=(F_1,\dots ,F_n)$ of $\Cn$ with a
finite number of zeros we define the Noether exponent $\nu(F)$.
We prove the Jacobi formula for all polynomials of degree
strictly less than $\sum_{i=1}^n (\deg F_i-1)-\nu(F)$.
\endabstract
\endtopmatter

\document
\head 		1. The Noether exponent
\endhead        
If $P=P(Z)$ is a complex polynomial in $n$ variables
$Z=(Z_1,\dots ,Z_n)$ then we denote by
$\tilde P=\tilde P(\tilde Z)$, $\tilde Z=(Z_0,Z)$ the
homogenization of $P$. If $\Ha$ is a set of homogeneous
polynomials in $n+1$ variables then we denote by $V(\Ha)$ the
subset of the complex projective space $\Pn$ defined by
equations $H=0$,\quad $H\in \Ha$.

The polynomial mapping $F=(F_1,\dots ,F_n):\Cn \to \Cn$
{\it has a finite number of zeros \/ }                  
if the set $V(\tilde F_1,\dots ,\tilde F_n)$ is finite.
We put $V_{\infty}(F)=V(\tilde F_1,\dots ,\tilde F_n,Z_0)$ and
call $V_{\infty}(F)$ the set of zeros of $F$ at infinity.
We identify $\Cn$ and $\Pn \setminus V(Z_0)$. Cleary $F$ has a
finite number of zeros if and only if the sets
$F^{-1}(0)\subset \Cn$ and $V_{\infty}(F)\subset \Pn$ are both finite.

\definition{Definition~(1.1)}
Let $F=(F_1,\dots ,F_n)$ be a polynomial mapping of $\Cn$ with a
finite number of zeros. By the Noether exponent of $F$ we mean
the smallest integer $\nu \geq 0$ such that the homogeneous
forms $Z_0^{\nu}$, $\tilde F_1$,\dots , $\tilde F_n$ satisfy
Noether's condition at every point of the set $V_{\infty}(F)$
(cf\.~Appendix).
\enddefinition
%
If $V_{\infty}(F)=\emptyset$ then $\nu (F)=0$.
If $V_{\infty}(F)\neq \emptyset$ and the hypersurfaces meet
transversally at any point of $V_{\infty}(F)$, then $\nu (F)=1$.
For any polynomial mapping $F=(F_1,\dots ,F_n)$ with a finite
number of zeros we put
$\mu (F)=\!\!\sum\limits_{z\in F^{-1}(0)}\!\!\mult_zF$
where $\mult_zF$ stands for the multiplicity of $F$ at $z$.
If $F^{-1}(0)=\emptyset$ then $\mu(F)=0$.

\medpagebreak
Let $d_i=\deg F_i$ for $i=1,\dots ,n$.

\proclaim{Proposition~(1.2)}
If $F$ has a finite number of zeros, then
$\nu(F)\leq \prod\limits_{i=1}^nd_i-\mu(F)$.
\endproclaim
%
\demo{Proof}
We have
$
\nu(F)\leq \max\{\,(\tilde F_1,\dots ,\tilde F_n)_p : p\in V_{\infty}(F)\,\}
$
(cf\.~Appendix (A5)). On the other hand, by Bezout's theorem
$\botsmash {
\sum\limits_{p\in V_{\infty}(F)}\!(\tilde F_1,\dots ,\tilde F_n)_p =
\prod\limits_{i=1}^nd_i-\mu(F) }
$
and (1.2) follows.
\enddemo

\remark{Remark(1.3)}
Let $k= \sharp V_{\infty}(F)$. Then a reasoning similar to the
above shows that $\nu(F)\leq \prod_{i=1}^nd_i-\mu(F)-k+1$.
\endremark

\proclaim{Proposition (1.4)}
Suppose that the polynomial mapping $F= (F_!,\dots ,F_n)$ has a
finite number of zeros and let $P$ be a polynomial belonging to
the ideal generated by~$F_1$,\dots , $F_n$ in the ring of
polynomials. Then there exist polynomials~$A_1$,\dots , $A_n$
such that
$P=A_1F_1+\dots +A_nF_n$ with $\deg A_iF_i\leq \deg P+\nu(F)$
for~$i=1$,\dots ,$n$.
\endproclaim

\demo{Proof}
The homogeneous forms
$Z_0^{\nu}\tilde P$, $\tilde F_1$,\dots , $\tilde F_n$ ($\nu=\nu(F)$)
satisfy Noether's conditions at
every point of $V(\tilde F_1,\dots ,\tilde F_n)$, then by
Noether's Fundamental Theorem (cf\. Appendix) there are
homogeneous forms~$\tilde A_1$,\dots , $\tilde A_n$ such that
$Z_0^{\nu}\tilde P=\tilde A_1\tilde F_1+\dots +\tilde A_n\tilde F_n$,
\quad
$\deg (\tilde A_i\tilde F_i)=\deg (Z_0^{\nu}\tilde P)=\nu +\deg P$.
We get (1.4) by putting $Z_0=1$.
\enddemo

For any $z=(z_1,\dots ,z_n)\in\Cn$ we put
$|z|=\max(|z_1|,\dots ,|z_n|)$. Recall that if $P:\Cn\to\Ce$ is
a polynomial of degree $d$ then there exist a constant $C>0$
such that $|P(z)|\leq C|z|^d$\quad for $|z|\geq 1$.

\proclaim{Proposition (1.5)}
Let $F=(F_1,\dots ,F_n)$ be a polynomial mapping with a finite
number of zeros. Then there exist positive constants $C$ and $R$
such that

$|F(z)|\geq C|z|^{\min(d_i)-\nu(F)}$\quad for $|z|\geq R$.
\endproclaim

\demo{Proof}
Since the fiber $F^{-1}(0)$ is finite then there are polynomials
$P_i(z_i)\not\equiv 0$ ($i=1$,\dots ,$n$) which belong to the ideal
generated by~$F_1$,\dots , $F_n$ in the ring of polynomials
(cf\.~\cite{4, p\.~23}). Let $m_i=\deg P_i(z_i)$. By (1.4) we can
write $P_i(z_i)=A_{i1}F_1+\dots +A_{in}F_n$,\quad
$\deg(A_{ij}F_j)\leq m_i+\nu(F)$. Hence there exist constants
$C>0$ and $R\geq 1$ such that
$$
|z_i|^{m_i}\leq C|z|^{m_i+\nu(F)-\min(d_i)}|F(z)|
$$
if $|z_i|\geq R$ for some $i\in \{\,1,\dots ,n\,\}$ and the
proposition follows.
\enddemo

\remark{Corollary~(1.6)}
If $\nu(F)<\min_{i=1}^n(d_i)$ then $F$ is proper i.e.
$\lim\limits_{|z|\to\infty}|F(z)| = {+}\infty$.
\endremark

To end with let us note two corollaries of propositions (1.2),
(1.4) and (1.5).

\remark{Corollary~(1.7)}
Let $F=(F_1,\dots ,F_n)$ be a polynomial mapping with a finite
number of zeros. Let $\mu=\mu(F)$. Then
\item{(1.7.1)}
(cf\. \cite 3, \cite{10}) there is a constant $C>0$ such that
$|F(z)|\geq C|z|^{\mu-\prod d_i+\min(d_i)}$ for large $|z|$.
\item{(1.7.2)} (cf\. \cite{11}) If $P$ belongs to the ideal
generated by $F_1$,\dots , $F_n$ in the ring of polynomials, then
$P=A_1F_1+\dots +A_nF_n$ with
$\deg(A_iF_i)\leq \prod_{i=1}^nd_i-\mu +\deg P$ for~$i=1$,\dots ,$n$.
\endremark

\head 		1. The Jacobi formula
\endhead        
Let $F=(F_1,\dots ,F_n)$ be a polynomial mapping such that the
fiber $F^{-1}(0)$ is finite and let $G:\Cn\to\Ce$ be a polynomial.
We denote by $\res_{F,z}(G)$ the residue of the meromorphic
differential form
$\dfrac{G(z)}{F_1(z)\dots F_n(z)}[dZ]$,\quad
$[dZ]=dZ_1\wedge\dots\wedge dZ_n$.

The definition and all properties of residues we need are given
in \cite{5}. Let us recall that if the Jacobian
$J_F=\det (\frac{\partial F_i}{\partial Z_j})$ is different from
zero at $z\in F^{-1}(0)$, then $\res_{F,z}(G)=\dfrac{G(z)}{J_F(z)}$.
The main result of this note is

\proclaim{Theorem~(2.1)}
Suppose that the polynomial mapping
$F=(F_0,\dots ,F_n):\Cn\to\Cn$ has a finite number of zeros.
Then the Jacobi formula
$$ \sum_{z\in F^{-1}(0)} \res_{F,z}(G)=0 \tag J $$
is satisfied for all polynomials $G:\Cn\to\Ce$ of degree
strictly less than $\sum_{i=1}^n(d_i-\nomathbreak 1)-\nu(F)$.
\endproclaim

Before giving the proof of (2.1) let us make some remarks.
If $F$ has no zeros at infinity i.e. if
$V_{\infty}(F)=\emptyset$  then $\nu(F)=0$ and (2.1) is reduced
to the Griffiths--Jacobi theorem (cf\. \cite{5}).
In~\cite{1}, \cite{2}, \cite{6} and~\cite{7} there are given
another generalisations of the Jacobi theorem. However, these
results does not imply ours.
If $\nu(F)\geq\sum_{i=1}^n(d_i-1)$ then the unique polynomial
satisfying the assumption of (2.1) is $G\equiv 0$.

\demo{Proof of ~(2.1)}
Let $\Omega$ be the meromorfic form in $\Pn$ given in $\Cn$ by
formula
$$ \Omega = \frac{G(z)}{F_1(z)\dots F_n(z)}[dZ]. $$
By Residue Theorem for $\Pn$ we get
$$  \sum_{z\in F^{-1}(0)}\!\res_{F,z}(G) =
   -\!\!\sum_{p\in V_{\infty}(F)}\!\res_p\Omega .  $$
It suffices to show, that $\res_p\Omega =0$ for all $p\in V_{\infty}(F)$.

Let $W=(W_1,\dots ,W_n)$ be an affine system of coordinates in
an affine neighbourhood of $p$ such that $W_1=0$ is the
hyperplane at infinity and $p$ has coordinates~$(0,\dots ,0)$.
Without loss of generality we may assume that
$Z_1=\frac{1}{W_1}$,
$Z_2=\frac{W_2+c_2}{W_1}$,\dots , $Z_n=\frac{W_n+c_n}{W_1}$.

Let
$P^{*}(W) =
W_1^dP(\frac{1}{W_1},\frac{W_2+c_2}{W_1},\dots ,\frac{W_n+c_n}{W_1})$
for any polynomial $P(Z)$ of degree~ $d$.
A simple calculation shows that near $p\in\Pn$:
$$
\Omega =\frac{-W_1^{\nu}G^{*}(W)}{F_1^{*}(W)\dots F_n^{*}(W)}[dW],\quad
\nu =\sum_{i=1}^n(d_i-1)-1-\deg G.
$$
By assumptions $\nu\geq\nu(F)$, therefore $W_1^{\nu}G^{*}(W)$
belongs to the local ideal generated
by~$F_1^{*}(W)$,\dots , $F_n^{*}(W)$. Consequently
$$
\res_p\Omega =
 \res_0\left(
  \frac{-W_1^{\nu}G^{*}(W)}{F_1^{*}(W)\dots F_n^{*}(W)}[dW]
\right) = 0
$$
and we are done.
\enddemo

\proclaim{Corollary~2.1}
If $F^{-1}(0)\neq\emptyset$ then
$\nu(F)\geq \sum_{i=1}^n(d_i-1)-\deg J_F$.
\endproclaim
\demo{Proof}
If $F^{-1}(0)\neq\emptyset$ then
$\sum\limits_{z\in F^{-1}(0)}\!\! \res_{F,z}(J_F)=\mu(F)\neq 0$,
consequently we cannot have $\deg J_F<\sum_{i=1}^n(d_i-1)-\nu(F)$.
\enddemo

\proclaim{Corollary~2.2}~(cf\. \cite{2})
If the hypersurfaces $\tilde F_i=0$ $(1\leq i\leq n)$ meet
transversally at infinity then (J) holds for any polynomial
of degree strictly less than~ $\sum_{i=1}^n(d_i-1)-1$.
\endproclaim
\demo{Proof}
If $\tilde F_i=0$ $(1\leq i\leq n)$ meet transversally
then $\nu(F)\leq 1$ and (2.2) follows immediately from (2.1).
\enddemo

\proclaim{Corollary~2.3}
Suppose that for any $p\in V_{\infty}(F)$: \par
(i) the hypersurfaces $\tilde F_i=0$  $(1\leq i\leq n)$
    have distinct tangent cones at $p$, \par
(ii) $\ord_p\tilde F_i\neq d_i$ \quad $(1\leq i\leq n)$. \newline
Then (J) holds for any polynomial of degree~${}\leq n-2$.
\endproclaim
\demo{Proof}
By (A6) we have
$\nu(F)\leq
 \max\{\,\sum\limits_{i=1}^n(\ord_p\tilde F_i-1)+1:p\in V_{\infty}(F)\,\}
 \leq \mathbreak
 \sum\limits_{i=1}^n(d_i-2)+1$
because $\ord_p(\tilde F_i)\leq d_i-1$ by (ii).
Consequently $\sum_{i=1}^n(d_i-1)-\nu(F)\leq n-1$ and it
suffices to use (2.1).
\enddemo

\example{Example}
Let $F(Z_1,Z_2)=(Z_1^{d_1}-1,Z_1Z_2+Z_2^{d_2})$ \quad
($d_1\geq 1$, $d_2\geq 2$).
Then condition (i) is satissfied but (ii) fails.
We have
$\botsmash{ \sum\limits_{z\in F^{-1}(0)}\res_{F,z}(1)=-1 }$,
hence condition (ii) is essential.
\endexample

\head 		Appendix: Noether's Conditions.
\endhead        
Let $H_0$, $H_1$,\dots , $H_n$ be homogeneous forms of $n+1$
variables such that the set $V=V(H_1,\dots ,H_n)$ is finite.
We denote by $\Hol_p$ the ring of holomorfic germs at~$p\in\Pn$.
Let $d_i=\deg H_i$ for $0\leq i\leq n$.

\proclaim{Max Noether's Fundamental Theorem}
The following two conditions are equivalent: \par
\item{(A1)} There is an equation  $H_0=A_1H_1+\dots A_nH_n$ (with $A_i$
            forms of degree $d_0-d_i$). \par
\item{(A2)} For any $p\in V$ there is an linear form $L$ such that
            $V\cap V(L)=\emptyset$ \newline
	    and
            $\dsize \frac{H_0}{L^{d_0}}\in \left(\frac{H_1}{L^{d_1}},\dots
             ,\frac{H_n}{L^{d_n}}\right)\Hol_p$.
\endproclaim

The proof of Noether's theorem follows easily
(cf\. \cite{4, p.~120}) from the affine version of the theorem
(cf\. \cite{12}) and from the following

\proclaim{(A3) Property}
If $H$ is a homogeneous form of $n+1$ variables such that
$V\cap V(H)=\emptyset$, then $H$ is not a zero-divisor modulo
ideal generated by $H_1$,\dots , $H_n$ in the ring of polynomials.
\endproclaim
\demo{Proof of (A3)}
If $V\cap V(H)=\emptyset$ then $H_1$,\dots , $H_n$, $H$ form the
sequence of parameters in the local ring $\Hol$ of holomorfic
functions at $0\in \Ce^{n+1}$, consequently $H$ is not a
zero-divisor ${}\mod(H_1,\dots ,H_n)\Hol$. Whence follows easily
(A3).
\enddemo

We say that the sequence $H_0$, $H_1$,\dots , $H_n$ satisfies
Noether's conditions at $p\in V$ if (A2) holds true. Let
$(H_1,\dots ,H_n)_p$ denotes the intersection number of~
$H_1$,\dots , $H_n$ at $p$ and let $\ord_pH$ be the order of $H$
at $p$. We have the following

\proclaim{Criteria for Noether's conditions}
The sequence $H_0$,\dots , $H_n$ satisfies \newline
Noether's conditions
at $p\in V$ if any of the following are true:

\item{(A4)} $H_1$,\dots , $H_n$ meet transversally at $p$ and
     $p\in V(H_0)$, \par
\item{(A5)} $\ord_pH_0 \geq (H_1,\dots ,H_n)_p$,\par
\item{(A6)} $H_1$,\dots , $H_n$ have distinct tangent cones at $p$ and
     $\ord_pH_0 \geq \sum_{i=1}^n (\ord_pH_i-\nomathbreak 1)+1$.
\endproclaim
\demo{Proof}
(A5) follows from the Mutiplicity theorem
(cf\.~\cite{8, p.~258}), (A6) is proved in \cite{9}
(Theorem~2.3),
(A4) is a special case both of (A5) and (A6).
\enddemo

\head 		Wyk\l{}adnik Noethera i formu\l{}a Jacobiego
\endhead        

\subhead Streszczenie
\endsubhead

{\eightpoint Dla ka{\accent95 z}dego odwzorowania wielomianowego
$F=(F_1,\dots ,F_n)$ przestrzeni $\Cn$ o sko\'nczonej liczbie zer
definiujemy wyk\l{}adnik Noethera $\nu(F)$
a nast\e{}pnie dowodzimy formu\l{}y Jacobiego dla wielomian\'ow
stopnia mniejszego od  $\sum_{i=1}^n (\deg F_i-1)-\nu(F)$. }
\enddocument